\documentclass[a4paper, 12pt]{amsart}
\usepackage{amsmath, amssymb, eucal, amscd, amstext, enumerate}
\usepackage{amsfonts,amsthm}
\usepackage{epsfig,graphicx,picinpar,subfigure,xypic}
\usepackage{pstricks}
\usepackage{fancyvrb}
\usepackage{fancyhdr}
\usepackage{mathrsfs}
\usepackage{indentfirst}
\usepackage{float}

\theoremstyle{plain}
\newtheorem{thm}{Theorem}[section]
\newtheorem{lem}[thm]{Lemma}

\newtheorem{prop}[thm]{Proposition}
\newtheorem{cor}[thm]{Corollary}
\newtheorem{defn}[thm]{Definition}
\newtheorem{rem}[thm]{Remark}
\newtheorem{rem-ntn}[thm]{Remark and Notation}
\newtheorem{ntn}[thm]{Notation}

\newenvironment{prf}{{\noindent \textbf{Proof:}\ }}{\hfill $\Box$\\ \smallskip}
\newenvironment{ack}{{\noindent \textbf{Acknowledgements} }}{\hfill \medskip}

\numberwithin{equation}{section}

\newcommand{\Cstar}{C^{\ast}}

\newcommand{\hull}{{\rm hull\ \!}}

\newcommand{\abs}[1]{\left\vert#1\right\vert}
\newcommand{\norm}[1]{\left\|#1\right\|}
\newcommand{\sph}{\mathfrak{S}_1}

\newcommand{\CP}{\mathcal{P}}
\newcommand{\CM}{\mathcal{M}}
\newcommand{\CN}{\mathcal{N}}
\newcommand{\CU}{\mathcal{U}}
\newcommand{\CV}{\mathcal{V}}
\newcommand{\CI}{\mathcal{I}}
\newcommand{\CJ}{\mathcal{J}}

\newcommand{\CL}{\mathcal{L}}
\newcommand{\CE}{\mathcal{E}}
\newcommand{\CH}{\mathcal{H}}

\newcommand{\CK}{\mathcal{K}}

\newcommand{\CA}{\mathcal{A}}
\newcommand{\CB}{\mathcal{B}}

\newcommand{\BC}{\mathbb{C}}

\newcommand{\BT}{\mathbb{T}}

\newcommand{\PH}{\mathcal{P}_{\mathcal{H}}}

\newcommand{\C}[1]{\mathcal{#1}}

\newcommand{\K}[1]{\emph{#1}}

\begin{document}

\title[On a geometric realization of $\Cstar$-algebras]{On a geometric realization of $\Cstar$-algebras}

\author{Xiao Chen}

\address[Xiao Chen]{Chern Institute of Mathematics, Nankai University, Tianjin 300071, China.}
\email{cxwhsdu@mail.nankai.edu.cn}

\thanks{}

\date{\today}



\begin{abstract}

Further to the functional representations of $\Cstar$-algebras proposed in \cite{CMP}, we
consider in this article the uniform K\"ahler bundle (in short, UKB) description of some
$\Cstar$-algebraic subjects. In particular, we obtain an one-to-one correspondence
between closed ideals of a $\Cstar$-algebra $\CA$ and full uniform K\"ahler subbundles over
open subsets of the base space of the UKB associated with $\CA$. In addition, we will
present a geometric description of the pure state space of hereditary $\Cstar$-subalgebras
and show that that if $\CB$ is a hereditary $\Cstar$-subalgebra of $\CA$, the UKB of $\CB$
is a kind of K\"ahler subbundle of the UKB of $\CA$. As a simple example, we consider
hereditary $\Cstar$-subalgebras of the $\Cstar$-algebra of compact operators on a Hilbert
space.
Finally, we remark that hereditary $\Cstar$-subalgebras also naturally can be
characterized as uniform holomorphic Hilbert subbundles.
\end{abstract}
\maketitle

\section{Introduction and preliminaries}

\medskip

For every commutative unital $\Cstar$-algebra $\CA$, it is well known that $\CA$ can
be faithfully represented as the algebra of the continuous functions on a compact
Hausdorff space. More precisely, Gelfand showed that $\CA$ can be viewed as the space
of continuous functions on the pure state space of $\CA$, equipped with the
$w^*$-topology. Gelfand's construction sets up an isomorphism between the
categories of commutative $\Cstar$-algebras (with unit) and compact Hausdorff spaces.
By means of this duality, any statement about commutative $\Cstar$-algebras can be
rephrased in the language of topological spaces, and vice versa.

It is rather natural to attempt a non-commutative version of this construction and try to
understand $\Cstar$-algebras in terms of continuous functions on certain spaces. A way
to obtain such a concrete realization had been proposed in \cite{CMP}. The space that
considered there is still the set of pure states as in the commutative case, equipped with
more complicated topological and geometric structures.
More precisely, the set $\CP$ will be viewed as a topological bundle of
infinite-dimensional K\"ahler manifolds over the base space of the
spectrum of the $\Cstar$-algebra (i.e. the topological space of unitary equivalence classes of non-zero
irreducible representations). In fact, each fiber of an irreducible representation consists of
pure states whose GNS representations are unitarily equivalent to that representation. In
this paper, we use this to give a geometric realization of some $\Cstar$-algebraic
subjects (such as closed ideals and hereditary $\Cstar$-subalgebras).

Firstly, we set up a one-to-one correspondence between closed ideals of
a general $\Cstar$-algebra $\CA$ and a class of uniform K\"ahler subbundles of pure
state space $\CP(\CA)$ of $\CA$. More precisely, we show that,
for any subbundle $(\CE,p_{\CA}|_{\CE},\C{X})$ of $(\CP(\CA),p_{\CA},\hat{\CA})$
such that $\C{X}$ is an open set of the spectrum $\hat{\CA}$ of $\CA$ and
$\CE=p_{\CA}^{-1}(\C{X})$, there is a close ideal $\CI$ in $\CA$ with
$\C{X}\cong\hat{\CA}^{\CI}$ and $\CP(\CI)\cong\CE$, and vice versa.

Secondly, we show that if $\CB$ is a hereditary $\Cstar$-subalgebra of a
$\Cstar$-algebra $\CA$, then $\CP(\CB)$ can be viewed as a kind of K\"ahler subbundle of $\CP(\CA)$
whose base space is a closed set of $\hat{\CA}$ and every fiber can be viewed as a closed K\"ahler submanifold of the corresponding fiber on
$\hat{\CA}$. Moreover, we discuss a kind of geometric structure on the pure state
space of a $\Cstar$-algebra and consider its relation to hereditary $\Cstar$-subalgebras.
Finally, as a simple example, we give a geometric characterization of hereditary
$\Cstar$-subalgebras of the $\Cstar$-algebra of compact operators on a Hilbert space.

We remark finally that, according to \cite{EK}, hereditary $\Cstar$-subalgebras
of a $\Cstar$-algebra $\CA$ also naturally can be characterized as another kind of
bundles associated to $\CA$, i.e., uniform holomorphic Hilbert bundles.

Before we start, let us first set some notations and recall from \cite{CMP} some results concerning
generalization of the Gelfand transform for non commutative unital $\Cstar$-algebras.

\begin{ntn}\label{ntn11}
All vector spaces in this article are over the complex field, unless stated
otherwise. We denote by $\langle\;|\;\rangle$ the inner product on a Hilbert space
$\CH$ (which is conjugate-linear in the first variable) and by $\CL(\CH)$ the space of all bounded linear operators on $\CH$.
\begin{itemize}
  \item We denote by $\CP_{\CH}$ the projective space of $\CH$.
        The element of $\CP_{\CH}$ is denoted by $[x]$, where $x\in\CH$ is a
        normalized representative of $[x]$. The unit sphere of $\CH$ is denoted by $\sph(\CH)$.
  \item We denote by $\CP(\CA)$ the set of all pure states on a $\Cstar$-algebra $\CA$
        endowed with $w^*$-topology and by $\hat{A}$ the spectrum of $\CA$ endowed with the Jacobson topology.
  \item We denote by
        $(\CP(\CA),p_{\CA},\hat{\CA})$ the uniform K\"ahler bundle (in short, UKB) associated with a $\Cstar$-algebra $\CA$,
        where $p_{\CA}:\CP(\CA)\longrightarrow\hat{\CA}$ is the natural projection given by the GNS representation (see Section 3 in \cite{CMP}).
\end{itemize}
\end{ntn}

Just as in the case of a commutative $\Cstar$-algebra, the function on the pure state
space $\CP(\CA)$ representing an element $a$ of $\CA$ will be the Gelfand
transform $f_a$ of $a$ given by
\begin{equation}\label{eq257}
f_a:\;\CP(\CA)\longrightarrow\BC,\quad \omega\longmapsto f_a(\omega):=\omega(a)\;.
\end{equation}

We introduce the main results in \cite{CMP} as follows.

Let $\C{C}^{\infty}(\CP)$ be the set of all smooth complex-valued functions on a
K\"ahler manifold $\CP$. If $\CP$ is the total space of a UKB, then
we still denote by $\C{C}^{\infty}(\CP)$ the set of all fiberwise-smooth complex-valued
functions on $\CP$ (see Section 3 in \cite{CMP} and Section 2.2 in \cite{KK}). Note that
$\C{C}^{\infty}(\CP)$ is a $*$-algebra with involution given by complex conjugation
$f\mapsto\bar{f}$.

If $\CA$ is a $\Cstar$-algebra, and $\CB=\CA\oplus\BC$ equipped with the canonical unital
$\Cstar$-algebraic structure, then $\hat{\CB}=\hat{\CA}\oplus\{\pi_0\}$,
$\CP(\CB)\cong\CP(\CA)\cup\{0\}\subseteq\CA^*$ and the fibre over $\pi_0\in\hat{\CB}$ is
$\{0\}$ (equipped with the trivial K\"ahler manifold structure). Moreover, for any
$f\in\C{C}^{\infty}(\CP(\CA))$, we extend $f$ to a function on $\CP(\CB)$ by setting
$f(0)=0$. Now, we have the following general form of \cite[Proposition 3.2]{CMP}.

\begin{thm}(\cite[Proposition 3.2]{CMP})\label{thm13}
Let $\CA$ be a $\Cstar$-algebra. Then
\begin{enumerate}
  \item The Gelfand transform $a\longmapsto f_a$ is a linear, involution
             preserving injection of $\CA$ into $\C{C}^{\infty}(\CP(\CA))$.
  \item The range of the Gelfand transform is the set denoted by
             $\C{K}_u(\CP(\CA))$ of $f\in\C{C}^{\infty}(\CP(\CA))$ such that $f$,
             $\bar{f}\star f$ and $f\star\bar{f}$ are uniformly continuous on
             $\CP(\CA)\cup\{0\}$ as well as $DDf=\bar{D}\bar{D}f=0$,
             where $D$ and $\bar{D}$ are respectively holomorphic and anti-holomorphic
             parts of covariant derivative of the K\"ahler metric defined on each fiber of $\CP(\CA)$, and the $\star$-product is as defined in \cite{CMP} (see also \cite[ Definition 2.1]{KK}).
  \item For any $a,b\in \CA$, one has $f_{ab}=f_a\star f_b$ and
             $\norm{a}^2=\sup_{\omega\in\CP(\CA)}(\bar{f_a}\star f_a)(\omega)$.
             By this norm, the $\C{K}_u(\CP(\CA))$ is a $\Cstar$-algebra which is $*$-isomorphic onto $\CA$.
\end{enumerate}
\end{thm}

\begin{defn}(\cite[Definition 2.1 and 2.2]{KK})\label{df13}
Two UKBs $(\CE,p,\C{X})$ and $(\CE',p',\C{X}')$ are \K{isomorphic} if there exists a pair
$(\psi,\phi)$ of homeomorphisms $\psi:\CE\longrightarrow\CE'$ and
$\phi:\C{X}\longrightarrow\C{X}'$, such that $p'\circ\psi=\phi\circ p$
and any restriction $\psi|_{p^{-1}(x)}:p^{-1}(x)\longrightarrow (p')^{-1}(\phi(x))$ is a
holomorphic K\"ahler isometry for any $x\in\C{X}$.
We call such a pair $(\psi,\phi)$ a \K{uniform K\"ahler isomorphism} between $(\CE,p,\C{X})$ and
$(\CE',p',\C{X}')$.
\end{defn}

\begin{cor}(\cite[Corollary 3.3]{CMP})\label{cor14}
If the UKBs $(\CP_1,p_1,\CB_1)$ and $(\CP_2,p_2,\CB_2)$ corresponding to the
$\Cstar$-algebras $\CA_1$ and $\CA_2$ are isomorphic, then the $\Cstar$-algebras $\CA_1$
and $\CA_2$ are themselves $*$-isomorphic.
\end{cor}

By the above results, we give a geometric structure on the pure state space and obtain a correspondence between algebra and geometry as follows:
\begin{equation*}
\xymatrix{
&{}\txt<8pc>{%
commutative $\Cstar$-algebra} \ar@2{<->}[r] \ar @{}[d] |{\bigcap}
&{}\txt<8pc>{%
locally compact Hausdorff space}
\ar @{}[d] |{\bigcap} \\
&{}\txt<8pc>{%
$\Cstar$-algebra} \ar@2{<->}[r]
&{}\txt<8pc>{%
UKB associated with a $\Cstar$-algebra}\,.
}
\end{equation*}

\bigskip

\section{Closed ideals, quotient algebras and hereditary subalgebras of $\Cstar$-algebras}

\medskip

In this section, we give a geometric characterization of quotient
$\Cstar$-algebras, closed ideals and hereditary subalgebras of a $\Cstar$-algebra.

First of all, we consider quotient $\Cstar$-algebras and closed ideals of a
$\Cstar$-algebra.
The results are quite obvious, but we include them here for completeness and comparison.

\begin{ntn}\label{ntn21}
Given a $\Cstar$-algebra $\CA$, let $\C{S}$ and $\CI$ be a subset and a
$\Cstar$-subalgebra of $\CA$ respectively.
\begin{itemize}
  \item Let
        $\hull(\C{S})$ be the set of primitive ideals of $\CA$ containing $\C{S}$.
  \item Let $\hat{\CA}^{\CI}$ be the set of
        $\pi\in\hat{\CA}$ such that $\pi(\CI)\neq0$.
  \item Let $\CP^{\CI}(\CA)$ be the set of pure states of
        $\CA$, which do not vanish on $I$.
\end{itemize}
\end{ntn}

\begin{defn}\cite[Definition 5.1]{Husemoller}\label{df22}
Let $(\CE,p,\C{X})$ be a bundle, and let $\CU$ be a subset of $\C{X}$. Then the \K{restriction of $(\CE,p,\C{X})$ to $\CU$}, denoted by $(\CE,p,\C{X})_{\CU}$, is the bundle $(\CE',p',\CU)$, where $\CE'=p^{-1}(\CU)$ and $p'=p|_{\CE'}$.
\end{defn}

\begin{ntn}\label{ntn23}
Let $\CA$ be a $\Cstar$-algebra. Set
\begin{itemize}
  \item $Bun_o(\CA):=\{(\CP(\CA),p_{\CA},\hat{\CA})_{\CU}\;|\;\CU\subseteq\hat{\CA}\text{ is an open subset}\}$.
  \item $Bun_c(\CA):=\{(\CP(\CA),p_{\CA},\hat{\CA})_{\CV}\;|\;\CV\subseteq\hat{\CA}\text{ is a closed subset}\}$.
\end{itemize}
\end{ntn}

\begin{rem-ntn}\label{rem-ntn24}
Note that $\CP_{\CH}$ is a \K{K\"ahler manifold}, and the K\"ahler distance
$d_{\CH}$ on $\CP_{\CH}$ is given by $d_{\CH}([x],[y])=\sqrt{2}\arccos\abs{\langle x|y\rangle}$ for any $[x],[y]\in\CP_{\CH}$ (see Appendix C in \cite{CMP}).

Besides, for each $[\pi]\in\hat{\CA}$, the fiber $p_{\CA}^{-1}([\pi])$ denoted
by $\CP_{\CA,[\pi]}$ is isomorphic, as a K\"ahler manifold, to the projective space of the Hilbert space $\CH_{\pi}$.
Indeed, the representation $\pi$ induces a K\"ahler isomorphism $\Phi_{\CA,[\pi]}: \CP_{\CA,[\pi]}\longrightarrow\CP_{\CH_{\pi}}$ defined by
$\Phi_{\CA,[\pi]}(\omega)=[x_{\omega}]\in\CP_{\CH_{\pi}}$ for any $\omega\in\CP_{\CA,[\pi]}$, where
$x_{\omega}$ is a canonical cyclic vector in $\CH_{\pi}$ satisfying
$\omega(a)=\langle x_{\omega}|\pi(a)x_{\omega} \rangle$ (see (3.1) in \cite{CMP} and Appendix D in \cite{CMP}). Moreover,
the \K{K\"ahler distance} in the $\CP(\CA)$ can be given as follows:
\begin{equation}\label{eq273}
d_{\CA}(\omega,\omega') = \left\{ \begin{array}{ll}
\sqrt{2}\arccos\abs{\langle x_{\omega}|x_{\omega'}\rangle}
 & \textrm{if $p_{\CA}(\omega)=p_{\CA}(\omega')$,}\\

3 & \textrm{otherwise,}
\end{array} \right.
\end{equation}
for any $\omega,\omega'\in\CP(\CA)$.
\end{rem-ntn}

\begin{prop}\label{prop25}
Let $\CA$ be a $\Cstar$-algebra. For any closed ideal $\CI$ of $\CA$, one has:

(a) The UKB $(\CP(\CI),p_{\CI},\hat{\CI})$ with respect to $\CI$ is uniformly
K\"ahler isomorphic to $(\CP(\CA),p_{\CA},\hat{\CA})_{\hat{\CA}^{\CI}}$.

(b) The quotient $\Cstar$-algebra $\CA/\CI$'s UKB
$(\CP(\CA/\CI),p_{\CA/\CI},\hat{\CA/\CI})$ is uniformly K\"ahler isomorphic to
$(\CP(\CA),p_{\CA},\hat{\CA})_{\hat{\CA}\setminus\hat{\CA}^{\CI}}$.
\end{prop}
\begin{prf}
We define by $f:\rho\mapsto\rho|_{\CI}$ the map from $\CP^{\CI}(\CA)$ to
$\CP(\CI)$, and by $g:\pi\mapsto\pi|_{\CI}$ the map
from $\hat{\CA}^{\CI}$ to $\hat{\CI}$. It is well known that $f$ and $g$ are both
homeomorphisms satisfying $g\circ p_{\CA}=p_{\CI}\circ f$ (see \cite[Section 2.11 and 3.2]{Dixmier}). For any
$[\pi]\in\hat{\CA}^{\CI}$ and $\rho\in\CP_{\CA,[\pi]}$, we have
$\CH_{\pi|_{\CI}}=\CH_{\pi}$ and $x_{\rho}=x_{\rho|_{\CI}}$. So
$\Phi_{\CA,[\pi]}\circ f\circ\Phi_{\CA,[\pi]}^{-1}$
is an identity map on $\CP_{\CH_{\pi}}$, namely, $f$ is a holomorphic K\"ahler
isometry. As a result, $(f,g)$ is a uniform K\"ahler isomorphism
between $(\CP(\CA),p_{\CA},\hat{\CA})_{\hat{\CA}^{\CI}}$ and $(\CP(\CI),p_{\CI},\hat{\CI})$. Hence (a) holds.

To prove (b), we denote by $h$ the quotient map from $\CA$ to $\CA/\CI$. It is well known that we can define by $f':\rho\mapsto\rho\circ h$ the homeomorphism from $\CP(\CA/\CI)$ to
$\CP(\CA)\setminus\CP^{\CI}(\CA)$ and by $g':\pi\mapsto\pi\circ h$ the homeomorphism
from $\hat{\CA/\CI}$ to $\hat{\CA}\setminus\hat{\CA}^{\CI}$ (see \cite[Section 2.11 and 3.2]{Dixmier}). Similar to the proof of (a), we can prove
that $(f',g')$ is a uniform K\"ahler isomorphism between
$(\CP(\CA),p_{\CA},\hat{\CA})_{\hat{\CA}\setminus\hat{\CA}^{\CI}}$ and
$(\CP(\CA/\CI),p_{\CA/\CI},\hat{\CA/\CI})$.
\end{prf}

\begin{prop}\label{prop26}
Let $\CA$ be a $\Cstar$-algebra. Then

(a) The correspondence
$\CI\longmapsto(\CP(\CA),p_{\CA},\hat{\CA})_{\hat{\CA}^{\CI}}$ is a bijection
from the set of closed ideals of $\CA$ onto $Bun_o(\CA)$, and $\CI=\ker \bigoplus_{[\pi]\in\hat{\CA}\setminus\hat{\CA}^{\CI}} \pi$.

(b) The correspondence
$\CA/\CI\longmapsto(\CP(\CA),p_{\CA},\hat{\CA})_{\hat{\CA}\setminus\hat{\CA}^{\CI}}$ is a bijection
from the set of the quotient $\Cstar$-algebras of $\CA$ onto $Bun_c(\CA)$.
\end{prop}
\begin{prf}
It follows from the fact that the correspondence $\CI\mapsto\hat{\CA}\setminus\hat{\CA}^{\CI}$ is
a bijection from the closed ideals of $\CA$ to the closed subsets of $\hat{\CA}$.
Moreover, we have
$\CI=\bigcap_{[\pi]\in\hat{\CA}\setminus\hat{\CA}^{\CI}}\ker\pi=\ker \bigoplus_{[\pi]\in\hat{\CA}\setminus\hat{\CA}^{\CI}} \pi$.
\end{prf}

In the rest of this section, we consider hereditary $\Cstar$-subalgebras.

Let $\CB$ be a nonzero hereditary $\Cstar$-subalgebra of a $\Cstar$-algebra
$\CA$. There is an injective map $\Xi$ from $\CP(\CB)$ to $\CP(\CA)$ such that
the image of $\tau$ in $\CP(\CB)$ under $\Xi$ is the unique extension of $\tau$ to $\CA$.
Moreover, every element in $\hat{\CB}$ has the form $[\pi_{\CB}]=[(\pi|_{\CB},\pi(\CB)\CH_{\pi})]$ for a unique $\pi\in\hat{\CA}^{\CB}$.

First, we discuss the relation between $\Xi(\CP(\CB))$ and $\CP^{\CB}(\CA)$.

\begin{prop}\label{prop27}
Let $\CB$ be a nonzero hereditary $\Cstar$-subalgebra of a $\Cstar$-algebra
$\CA$ and $[\pi]\in\hat{\CA}^{\CB}$. Then

(a) if $\pi(\CB)\CH_{\pi}=\CH_{\pi}$, there is a unique bijection
             $\Theta: \CP_{\CA,[\pi]}\rightarrow
             \CP_{\CB,[\pi_{\CB}]}$ such that $\Theta(\rho)=\rho|_{\CB}$
             and $\Theta=\Xi^{-1}$.

(b) if $\pi(\CB)\CH_{\pi}\subsetneq\CH_{\pi}$, there is a unique surjection
             $\Theta: \CP_{\CA,[\pi]}\cap\CP^{\CB}(\CA)\rightarrow(0,1]\times
              \CP_{\CB,[\pi_{\CB}]}$ such that $\rho|_{\CB}=t\rho'$ if $\Theta(\rho)=(t,\rho')$.

\end{prop}
\begin{prf}
(a) If $\pi(\CB)\CH_{\pi}=\CH_{\pi}$, then $\CH_{\pi_{\CB}}=\CH_{\pi}$. So,
for any $\rho\in\CP_{\CA,[\pi]}$, we have $\rho|_{\CB}\in\CP_{\CB,[\pi_{\CB}]}$ and
$\Xi(\Theta(\rho))=\rho$. On the other hand, for each
$\rho'\in\CP_{\CB,[\pi_{\CB}]}$, we have $\Theta(\Xi(\rho'))=\rho'$.
Hence (a) holds.

(b) We first show that $\Theta$ is well-defined. For any $\rho\in\CP_{\CA,[\pi]}\cap\CP^{\CB}(\CA)$, by \cite[Corollary 5.5.3]{Murphy}, there exists a pair $(t,\rho')\in(0,1]\times\CP(\CB)$ such that
$\rho|_{\CB}=t\rho'$. If there is another pair
$(s,\omega)\in(0,1]\times\CP_{\CB,[\pi_{\CB}]}$ such that
$\rho|_{\CB}=s\omega$, then $t=s$ and $\rho'=\omega=\frac{1}{t}\rho|_{\CB}$ as $\norm{t\rho'}=\norm{s\omega}$ as well as
$\norm{\rho'}=\norm{\omega}=1$. Next, we show that $\Theta$ is a surjection.
For any $(t,\rho')\in(0,1)\times\CP_{\CB,[\pi_{\CB}]}$, let
$(\CH_{\rho'},\pi_{\rho'})$ be the GNS representation of $\CB$ associated with
$\rho'$. Then $\pi_{\CB}$ is unitarily equivalent to
$(\CH_{\rho'},\pi_{\rho'})$. Let $x\in\sph(\pi(\CB)\CH_{\pi})$ be the canonical
cyclic vector with $\rho'(b)=\langle x|\pi_{\CB}(b)x\rangle$. Set
$y=\sqrt{t}x\in\pi(\CB)\CH_{\pi}$, so that $\norm{y}=\sqrt{t}$. Since $\pi(\CB)\CH_{\pi}\subsetneq\CH_{\pi}$,
 we can choose a $z\in(\pi(\CB)\CH_{\pi})^{\bot}$ satisfying
$\norm{z}=\sqrt{1-t}$. Let $P$ be a projection $\CH_{\pi}\mapsto\pi(\CB)\CH_{\pi}$ and
put $h=y+z\in\CH_{\pi}$. It's clear that $\norm{h}=1$ and $P(h)=y$. For any
$v,w\in\CH_{\pi}$ and $b\in\CB$,
\begin{eqnarray*}
  \langle w|\pi(b)P(v)\rangle &=& \langle\pi(b^*)w|P(v)\rangle
        =\langle P\pi(b^*)w|v\rangle=\langle\pi(b^*)w|v\rangle \\
  {} &=& \langle w|\pi(b)(v)\rangle=\langle w|P\pi(b)(v)\rangle .
\end{eqnarray*}
So $P\in\pi(\CB)'$.
Let $\rho(a)$ be $\langle h|\varphi(a)h\rangle\;(a\in\CA)$. Then, for any
$b\in\CB$,
\begin{eqnarray*}
  \rho(b) &=& \langle h|\pi(b)h\rangle=\langle h|P\pi(b)h\rangle
        =\langle P(h)|\pi(b)P(h)\rangle=\langle y|\pi(b)y\rangle \\
  {} &=& \langle y|\pi_{\CB}(b)y\rangle
        =\langle\sqrt{t}x|\pi_{\CB}(b)(\sqrt{t}x)\rangle=t\rho'(b),
\end{eqnarray*}
which implies that $\rho|_{\CB}=t\rho'$. Besides, it is obvious that $\rho$
belongs to $\CP_{\CA,[\pi]}\cap\CP^{\CB}(\CA)$. Based on the discussion above, we
know that $\rho$ belongs to pre-image $\Theta^{-1}(t,\rho')$. So (b) holds.
\end{prf}

\begin{rem}\label{rem28}
If $\pi\in\hat{\CA}^{\CB}$ and
$\Delta:=\Theta^{-1}(\{1\}\times \CP_{\CB,[\pi_{\CB}]})$, then
$\Theta|_{\Delta}$ can be identified with the inverse map of $\Xi$.
\end{rem}

Moreover, we can describe the relation between $\Xi(\CP(\CB))$ and $\CP^{\CB}(\CA)$ by endowing a kind of geometric structure on their pure state spaces.

For any $\mu\in\CP(\CA)$ and $t\in(0,+\infty)$, set
$\BT:=\{\lambda\in\BC\;|\;\abs{\lambda}=1\}$,
$S(\mu;t):=\{\nu\in\CP(\CA)\;|\;d_{\CA}(\mu,\nu)=t\}$,
$D(\mu;t):=\{\nu\in\CP(\CA)\;|\;d_{\CA}(\mu,\nu)<t\}$ and
$B(\mu;t):=\{\nu\in\CP(\CA)\;|\;d_{\CA}(\mu,\nu)\leqslant t\}$.

\begin{thm}\label{thm29}
Let $\CB$ be a nonzero hereditary $\Cstar$-subalgebra of a $\Cstar$-algebra
$\CA$, and set $\kappa:=\frac{\sqrt{2}\pi}{2}$. For any
$[\pi]\in\hat{\CA}^{\CB}$, $\mu\in\CP_{\CB,[\pi_{\CB}]}$ and
$t\in(0,\kappa)$, we have

(a) $$p_{\CA}^{-1}(\hat{\CA}\setminus\hat{\CA}^{\CB})=\{\rho\in\CP(\CA)\;|\;d_{\CA}(\rho,\Xi(\CP(\CB)))=3\},$$
$$\CP^{\CB}(\CA)=\{\rho\in\CP(\CA)\;|\;d_{\CA}(\rho,\Xi(\CP(\CB)))<\kappa\},$$
and
$$p_{\CA}^{-1}(\hat{\CA}^{\CB})\setminus\CP^{\CB}(\CA)
=\{\rho\in\CP(\CA)\;|\;d_{\CA}(\rho,\Xi(\CP(\CB)))=\kappa\}.$$

(b)
$\CP_{\CA,[\pi]}=B(\Xi(\mu);\kappa)$,
$\CP_{\CA,[\pi]}\cap\CP^{\CB}(\CA)= D(\Xi(\mu);\kappa)$, and
             $\CP_{\CA,[\pi]}\setminus\CP^{\CB}(\CA)= S(\Xi(\mu);\kappa)$.

(c) $\Theta^{-1}(\cos^2(\frac{t}{\sqrt{2}}),\mu)=S(\Xi(\mu);t)$
and there is a bijection
$$\Upsilon:S(\Xi(\mu);t)\rightarrow\BT\times(\CP_{\CA,[\pi]}\setminus\CP^{\CB}(\CA)).$$

\end{thm}
\begin{prf}
According to Proposition~\ref{prop27}, we obtain two surjections
$$\hat{\Theta}:\CP(\CA)^{\CB}\rightarrow\CP(\CB), \rho\mapsto\rho',$$ and
$$t_{\CB}:\CP(\CA)^{\CB}\rightarrow(0,1], \rho\mapsto t,$$ where $(t,\rho')=\Theta(\rho)$. It is apparent that
$t_{\CB}$ can be extended to $\CP(\CA)$ by setting $$t_{\CB}(\rho)=0, \forall\rho\in\CP(\CA)\setminus\CP(\CA)^{\CB}.$$

For any $\rho\in\CP(\CA)^{\CB}$, let
$(\CH_{\rho},\pi_{\rho},x_{\rho})$ be the GNS representation of $\CA$ associated
to $\rho$. One may identify
$\CH_{\hat{\Theta}(\rho)}=\pi_{\rho}(\CB)\CH_{\rho}\subseteq\CH_{\rho}$
and $\pi_{\hat{\Theta}(\rho)}(b)=\pi_{\rho}(b)|_{\pi_{\rho}(\CB)\CH_{\rho}}$.
Moreover, there is a unique
$w_{\rho}\in\sph(\CH_{\rho}\ominus\CH_{\hat{\Theta}(\rho)})$ such that
$$x_{\rho}=\sqrt{t_{\CB}(\rho)}x_{\hat{\Theta}(\rho)}+\sqrt{1-t_{\CB}(\rho)}w_{\rho},$$
which implies that
$$d_{\CA}(\rho,\Xi(\hat{\Theta}(\rho)))=\sqrt{2}\arccos\sqrt{t_{\CB}(\rho)}.$$
Consequently, $$\CP(\CA)^{\CB}=\bigcup_{\omega\in\CP(\CB)}D(\Xi(\omega);\kappa)\subseteq p_{\CA}^{-1}(\hat{\CA}^{\CB})=\bigcup_{\omega\in\CP(\CB)}B(\Xi(\omega);\kappa).$$
From the above argument, it can be easily checked that (a) holds.

For any $\mu\in\CP_{\CB,[\pi_{\CB}]}$ and $\rho\in\CP_{\CA,[\pi]}$,
one may identify $\CH_{\rho}=\CH_{\pi}$ and $\CH_{\mu}=\CH_{\pi_{\CB}}$. So
there is a unique $w_{\rho}\in\sph(\CH_{\pi}\ominus\CH_{\pi_{\CB}})$ such that
$x_{\rho}=\sqrt{t_{\CB}(\rho)}x_{\mu}+\sqrt{1-t_{\CB}(\rho)}w_{\rho}$.
For any $t\in(0,\kappa)$, we have
$\Theta^{-1}(\cos^2(\frac{t}{\sqrt{2}}),\mu)$ is the subset of $\CP_{\CA,[\pi]}$,
which is as follow:
$$\{\rho\in\CP_{\CA,[\pi]}\;|\;x_{\rho}=
\cos(\frac{t}{\sqrt{2}})x_{\mu}+\sqrt{1-\cos^2(\frac{t}{\sqrt{2}})}w, \forall w\in\sph(\CH_{\pi}\ominus\CH_{\pi_{\CB}})\},$$
which implies that
$$\Theta^{-1}(\cos^2(\frac{t}{\sqrt{2}}),\mu)\cong\sph(\CH_{\pi}\ominus\CH_{\pi_{\CB}})
\cong\BT\times\CP_{\CH_{\pi}\ominus\CH_{\pi_{\CB}}}.$$ In addition, we have that $$d_{\CA}(\rho,\Xi(\mu))=\sqrt{2}\arccos\sqrt{t_{\CB}(\rho)},
\Phi_{\CA,[\pi]}(\CP_{\CA,[\pi]}\setminus\CP^{\CB}(\CA))=
\CP_{\CH_{\pi}\ominus\CH_{\pi_{\CB}}}.$$ So (b) and (c) hold.
\end{prf}

Next, we characterize the relationship between $\CP(\CB)$ and $\CP(\CA)$ by using
the notion of K\"ahler bundle. Before that, we recall the following well-known result.

\begin{lem}\label{lem210}
Let $\CI$ be a closed ideal which is generated by a nonzero hereditary $\Cstar$-subalgebra
$\CB$ of a $\Cstar$-algebra $\CA$. Then $\hat{\CA}^{\CB}=\hat{\CA}^{\CI}$ is an
open set in $\hat{\CA}$. Moreover, one has $\CI=\bigcap_{\CJ\in\hull(\CB)}\CJ$.
\end{lem}

\begin{defn}(\cite[Definition 8.5]{Upmeier})\label{df211}
A closed (or open) subset $\CN$ of a K\"ahler Hilbert manifold $\CM$ is called a
\K{closed (or open) submanifold} if for any $x\in\CN$, there exists a chart
$(\CV,b_{\CV},{\CH}_{\CV})$ of $x$ and a closed subspace
$\CE$ of $\CH_{\CV}$ such that $b_{\CV}(\CN\cap\CV)=\CE\cap b_{\CV}(\CV)$.
\end{defn}

\begin{lem}\label{lem212}
If $\CH$ is a Hilbert space and $\CM$ is a closed subspace of $\CH$, then
$\CP_{\CM}$ is a closed K\"ahler submanifold of $\CP_{\CH}$.
\end{lem}
\begin{prf}
For any $[\xi]\in\CP_{\CM}$, let $(\CV_{\xi},b_{\xi},\CH_{\xi})$ be a canonical
chart of $[\xi]$ (see Appendix C in \cite{CMP}).
Since $b_{\xi}(\CP_{\CM}\cap\CV_{\xi})=\CM\cap\CH_{\xi}$, it follows that
$\CP_{\CM}$ is a K\"ahler submanifold of $\CP_{\CH}$.
Moreover, it can be easily checked that $\CP_{\CM}$ is closed under the
K\"ahler topology induced by the K\"ahler distance $d_{\CH}$.
\end{prf}

\begin{thm}\label{thm213}
Let $\CB$ be a nonzero hereditary $\Cstar$-subalgebra of a $\Cstar$-algebra $\CA$. Then

(a) $(\CP(\CA),p_{\CA},\hat{\CA})^{\CB}:=(\Xi(\CP(\CB)),p_{\CA},\hat{\CA}^{\CB})$
is a K\"ahler subbundle of $(\CP(\CA),p_{\CA},\hat{\CA})$ such that
$\Xi(\CP_{\CB,[\pi_{\CB}]})$ is a closed K\"ahler submanifold of
$\CP_{\CA,[\pi]}$ for all $[\pi]\in\hat{\CA}^{\CB}$. Moreover,
$(\CP(\CB),p_{\CB},\hat{\CB})$ is uniformly K\"ahler isomorphic to $(\CP(\CA),p_{\CA},\hat{\CA})^{\CB}$.

(b) $\CB$ is a closed ideal if only and if
$(\CP(\CA),p_{\CA},\hat{\CA})^{\CB}=(\CP(\CA),p_{\CA},\hat{\CA})_{\hat{\CA}^{\CB}}$.
\end{thm}
\begin{prf}
(a) We denote by $\Psi$ the canonical homeomorphism from $\hat{\CA}^{\CB}$ to
$\hat{\CB}$. Let $\Xi$ be a homeomorphism from $\CP(\CB)$ to $\Xi(\CP(\CA))$
satisfying $\Psi^{-1}\circ p_{\CB}=p_{\CA}\circ\Xi$. By Proposition~\ref{prop27},
for any $[\pi]\in\hat{\CA}^{\CB}$, we have
$\Xi(\CP_{\CB,[\pi_{\CB}]})\subseteq\CP_{\CA,[\pi]}$.
So $(\CP(\CA),p_{\CA},\hat{\CA})^{\CB}$ is a subbundle of
$(\CP(\CA),p_{\CA},\hat{\CA})$. Since
$\Phi_{\CA,[\pi]}(\CP_{\CA,[\pi]})=\CP_{\CH_{\pi}}$ and
$\Phi_{\CB,[\pi_{\CB}]}(\CP_{\CB,[\pi_{\CB}]})=\CP_{\CH_{\pi_{\CB}}}$, according to Lemma~\ref{lem212}, we know that
$\Xi(\CP_{\CB,[\pi_{\CB}]})$ is a closed
K\"ahler submanifold of $\CP_{\CA,[\pi]}$.
Moreover, according to the above argument, it is obvious that $(\Xi,\Phi^{-1})$ is a
uniform K\"ahler isomorphism between $(\CP(\CB),p_{\CB},\hat{\CB})$ and
$(\CP(\CA),p_{\CA},\hat{\CA})^{\CB}$. Hence (a) holds.

(b) It is obvious by Proposition~\ref{prop25}(a) and Lemma~\ref{lem212}.
\end{prf}

Finally, we give a concrete example: $\CA$ is $\CK(\CH)$ consisting of all compact operators on a Hilbert space $\CH$.

It can be easily proved that the UKB associated with $\CA$ has only one fibre denoted by $\CP$ because $\hat{\CA}=\{[(\CH,i)]\}$ (see \cite[Example 5.1.1]{Murphy}), where $i$ is the
inclusion map from $\CA$ to $\CL(\CH)$. We denote by $\Phi$ the canonical K\"ahler
isomorphism from $\CP$ to $\CP_{\CH}$ and by $\pi$ the quotient map from $\CH$ to
$\PH$. For any $[\xi]\in\PH$, set
$\{\xi\}^{\bot}:=\{x\in\CH\;|\;\langle\xi|x\rangle=0\}$ and denote by $[S]$ the closed
linear span of a set $S\subset\CH$.

Note that $\{\xi\}^{\bot}\cong T_{[\xi]}\Phi(\CP)$,
where $T_{[\xi]}\Phi(\CP)$ is the \K{holomorphic tangent space} (See Appendix A in \cite{CMP} and Section 4 in \cite{Upmeier}) at $[\xi]$ in the K\"ahler manifold $\Phi(\CP)$.

\begin{prop}\label{prop214}
Let $\CA$ be $\CK(\CH)$. Then $\CB\subseteq\CA$ is a hereditary $\Cstar$-subalgebra if and only
if the UKB associated with $\CB$ has only one fibre $\CP'$ which is a closed K\"ahler
submanifold of $\CP$ satisfying $$\pi([\sum_{x\in\CP'} T_{\Phi(x)}\Phi(\CP')])\subset\Phi(\CP').$$ Moreover, one has $\CP'=\Phi^{-1}(\CP_{\CM})$ and $\CB=\CK(\CM)$, where $\CM:=[\sum_{x\in\CP'} T_{\Phi(x)}\Phi(\CP')]$.
\end{prop}

\begin{prf}
It is well known that $\CB$ is a hereditary $\Cstar$-subalgebra of $\CK(\CH)$ if
and only if there is a closed subspace $\CM$ of $\CH$ such that $\CB=\CK(\CM)$.
So the forward implication is obvious considering $\hat{\CA}=\{[(\CH,i)]\}$.

Conversely, assume that $\CP'$ is a closed K\"ahler submanifold of $\CP$ such that
$$\pi([\sum_{x\in\CP'} T_{\Phi(x)}\Phi(\CP')])\subset\Phi(\CP').$$
Choose an arbitrary point
$[\xi]$ in $\Phi(\CP')\subset\Phi(\CP)$, and let $(\CV_{\xi},b_{\xi},\CH_{\xi})$ be a chart of $[\xi]$.
Since $\CP'$ is a closed K\"ahler submanifold of $\CP$, there is a chart
$(\CV_{\xi}\cap\Phi(\CP'),b_{\xi}|_{\CV_{\xi}\cap\Phi(\CP')},\C{W}_{\xi})$ of $\Phi(\CP')$ such that
$b_{\xi}(\CV_{\xi}\cap\Phi(\CP'))=\C{W}_{\xi}\cap b_{\xi}(\CV_{\xi})$.
By the assumption, we have $$\pi(T_{[\eta]}\Phi(\CP'))\subset\Phi(\CP'),\forall [\eta]\in\Phi(\CP').$$
It follows from the definition of $b_{\xi}$ that
$$\pi(\{\eta\}\oplus T_{[\eta]}\Phi(\CP'))\subset\Phi(\CP'),\forall [\eta]\in\Phi(\CP').$$
For each $[\zeta]\in\Phi(\CP')$, there exists $[\eta]\in\Phi(\CP')$ such that $\zeta\in T_{[\eta]}(\Phi(\CP'))$. Hence
$$\Phi(\CP')\subset\bigcup_{[\eta]\in\Phi(\CP')}\pi(T_{[\eta]}\Phi(\CP')).$$
Since $$\bigcup_{[\eta]\in\Phi(\CP')}\pi(T_{[\eta]}\Phi(\CP'))\subset
\pi([\sum_{x\in\CP'} T_{\Phi(x)}\Phi(\CP')]),$$ we have $$\Phi(\CP')\subset\pi([\sum_{x\in\CP'} T_{\Phi(x)}\Phi(\CP')]).$$
Consequently, $$\Phi(\CP')=\pi([\sum_{x\in\CP'} T_{\Phi(x)}\Phi(\CP')]).$$ Thus, we obtain that $\CP'=\Phi^{-1}(\CP_{\CM})$ and $\CB=\CK(\CM)$, where $\CM:=[\sum_{x\in\CP'} T_{\Phi(x)}\Phi(\CP')]$.
\end{prf}

\begin{rem}\label{rem215}
Let $\CA$ be a $\Cstar$-algebra. Consider $\CA$ as a right Hilbert module over
itself, i.e., with the $\CA$-valued inner product $(a,b)\mapsto a^{*}b$ $(\forall
a,b\in\CA)$. According to \cite{EK}, $\CA$ can be view as a uniform holomorphic
Hilbert bundle $\{\CH(\CA)_{\rho}\}_{\rho\in\CP(\CA)\cup\{0\}}$, where
$\CH(\CA)_{\rho}$ is the GNS Hilbert space associated to $\rho$. Hence, each
$\Cstar$-algebra $\CA$ can induce two bundles
$\{\CH(\CA)_{\rho}\}_{\rho\in\CP(\CA)\cup\{0\}}$ and $(\CP(\CA),p_{\CA},\hat{\CA})$
(the base space of the former bundle is the total space of the latter bundle).

It is well known that the correspondence $\CL\mapsto \CL\cap\CL^{*}$ is a
bijection from the set of closed left ideals of $\CA$ onto the set of hereditary
$\Cstar$-subalgebras of $\CA$. We denote by $\CL(\CB)$ the closed left ideal
associated to a hereditary $\Cstar$-subalgebra $\CB$ of $\CA$. Since
$\CL^{*}$ is a closed right ideal in $\CA$, $\CL^{*}$ can be considered as a right
Hilbert sub-$\CA$-module of $\CA$. So it follows from \cite{EK} that $\CL$ can
induce a subbundle $\{\CH(\CL^{*})_{\rho}\}_{\rho\in\CP(\CA)\cup\{0\}}$ of
$\{\CH(\CA)_{\rho}\}_{\rho\in\CP(\CA)\cup\{0\}}$. Consequently, $\CB$ corresponds to such subbundle.

In brief, each hereditary $\Cstar$-subalgebra $\CB$ of $\CA$ can be described
as two different bundles associated to $\CA$ respectively, i.e., a uniform K\"ahler subbundle
$(\CP(\CA),p_{\CA},\hat{\CA})^{\CB}$ and a uniform holomorphic Hilbert subbundle
$\{\CH(\CL(\CB)^{*})_{\rho}\}_{\rho\in\CP(\CA)\cup\{0\}}$.
\end{rem}

\bigskip

\begin{ack}
The author would like to sincerely appreciate his supervisor Professor
Chi-Keung Ng (Nankai University) for his guidance and valuable suggestions. The
author also wants to thank the referees for valuable comments which improve this
paper a lot. Meanwhile, the author would like to thank Mr.Jianze Li (Nankai University), Mr.Qun Zhang (University of South Australia), Mr.Xiaoyu Li (Harbin Institute of Technology) and Ms.Guoping Zhao (Zhejiang University) for their generous help.
\end{ack}

\end{document}